\documentclass[a4paper,12pt]{article}
\usepackage{amsthm,amsmath,amssymb,latexsym,bm}

\newtheorem{theorem}{Theorem}[section]
\newtheorem{corollary}[theorem]{Corollary}
\newtheorem{proposition}[theorem]{Proposition}
\newtheorem{lemma}[theorem]{Lemma}
\newtheorem{remark}[theorem]{Remark}
\numberwithin{equation}{section}
\title{A $q$-analogue of the Drinfeld-Sokolov hierarchy of type $A$ and $q$-Painlev\'{e} system}
\author{Takao Suzuki \thanks{Department of Mathematics, Kinki University, 3-4-1, Kowakae, Higashi-Osaka, Osaka 577-8502, Japan. E-mail: suzuki@math.kindai.ac.jp}}
\date{}

\begin{document}

\maketitle

\begin{abstract}
In this article, we propose a $q$-analogue of the Drinfeld-Sokolov (DS) hierarchy of type $A$.
We also discuss its relationship with the $q$-Painlev\'{e} VI equation and the $q$-hypergeometric function.


\end{abstract}

\section{Introduction}

A relationship between Painlev\'{e} systems and infinite-dimensional integrable hierarchies has been studied; e.g. \cite{FN,FS1,KK2,NY1,T2}.
In a recent work \cite{FS2,S2}, we derive a coupled Painlev\'{e} VI system, which admits a particular solution in terms of the hypergeometric function ${}_nF_{n-1}$ \cite{S1,T3}, from the DS hierarchy of type $A$.
In this article, we give a $q$-difference analogue of the above result.
Namely, we propose a $q$-analogue of the DS hierarchy of type $A$ and derive a $q$-Painlev\'{e} system which admits a $q$-hypergeometric solution.

The DS hierarchies are extensions of the KdV hierarchy for affine Lie algebras \cite{DS,FHM,GHM}.
For type $A_{N-1}^{(1)}$ among them, the hierarchies are characterized by partitions of a natural number $N$.
In this article, we propose a $q$-DS hierarchy corresponding to the partition $(n,\ldots,n)$ of $mn$.
Note that its formulation is based on two preceding works \cite{KNY2}, equivalent to the case $m=1$, and \cite{KK1}, equivalent to the case $n=1$.

The coupled Painlev\'{e} VI system given in \cite{S2} (or \cite{T4}) is derived from the DS hierarchy corresponding to the partition $(n,n)$ by a similarity reduction.
Therefore, in this article, we investigate a similarity reduction of the $q$-DS hierarchy corresponding to the partition $(n,n)$; we denote it by $q$-$P_{(n,n)}$.
The system $q$-$P_{(2,2)}$ coincides with the $q$-Painlev\'{e} VI equation given in \cite{JS}.
Thus we can regard $q$-$P_{(n,n)}$ as a higher order generalization of the $q$-Painlev\'{e} VI equation.
The system $q$-$P_{(n,n)}$ also admits a particular solution in terms of the $q$-hypergeometric function ${}_n\phi_{n-1}$.

The $q$-Painlev\'{e} VI equation is a system of $q$-difference equations
\begin{equation}\begin{split}\label{Eq:qP6_JS}
	\frac{x(t)x(q^{-1}t)}{\alpha_3\alpha_4} &= \frac{(y(q^{-1}t)-t\beta_1)(y(q^{-1}t)-t\beta_2)}{(y(q^{-1}t)-\beta_3)(y(q^{-1}t)-\beta_4)},\\
	\frac{y(t)y(q^{-1}t)}{\beta_3\beta_4} &= \frac{(x(t)-t\alpha_1)(x(t)-t\alpha_2)}{(x(t)-\alpha_3)(x(t)-\alpha_4)},
\end{split}\end{equation}
where the parameters satisfy
\begin{equation*}
	\frac{\beta_1\beta_2}{\beta_3\beta_4} = q^{-1}\frac{\alpha_1\alpha_2}{\alpha_3\alpha_4}.
\end{equation*}
It is given as the compatibility condition of a system of linear $q$-difference equations
\begin{equation}\begin{split}\label{Eq:Lax_JS}
	&Y(q^{-1}z,t) = \mathcal{A}(z,t)Y(z,t),\quad \mathcal{A}(z,t) = \mathcal{A}_0(t) + \mathcal{A}_1(t)z + \mathcal{A}_2(t)z^2,\\
	&Y(z,q^{-1}t) = \frac{z(zI+\mathcal{B}_0(t))}{(z-q^{-1}t\alpha_1)(z-q^{-1}t\alpha_2)}Y(z,t).
\end{split}\end{equation}
The coefficient matrices $\mathcal{A}_i(t)$ $(i=0,1,2)$ satisfy
\begin{equation*}\begin{split}
	&\mathcal{A}_2(t) = \begin{bmatrix}\kappa_1&0\\0&\kappa_2\end{bmatrix},\quad \text{$\mathcal{A}_0(t)$ has eigenvalues $t\theta_1,t\theta_2$},\\
	&\det\mathcal{A}(z,t) = \kappa_1\kappa_2(z-t\alpha_1)(z-t\alpha_2)(z-\alpha_3)(z-\alpha_4),
\end{split}\end{equation*}
where
\begin{equation*}
	\kappa_1 = \frac{q}{\beta_3},\quad \kappa_2 = \frac{1}{\beta_4},\quad \theta_1 = \frac{\alpha_1\alpha_2}{\beta_1},\quad \theta_2 = \frac{\alpha_1\alpha_2}{\beta_2}.
\end{equation*}
The dependent variables $x(t)$ and $y(t)$ are given by
\begin{equation*}
	x(t) = -\frac{(\mathcal{A}_0(t))_{12}}{(\mathcal{A}_1(t))_{12}},\quad y(t) = \frac{q(x(t)-t\alpha_1)(x(t)-t\alpha_2)}{(\mathcal{A}(x(t),t))_{11}}.
\end{equation*}
We derive the system \eqref{Eq:qP6_JS} from $q$-$P_{(2,2)}$ with the aid of a $q$-Laplace transformation for a system of linear $q$-difference equations.

\begin{remark}
The $q$-Painlev\'{e} VI equation \eqref{Eq:qP6_JS} was derived from $q$-$\mathfrak{gl}_3$ hierarchy via the $q$-Laplace transformation in {\rm\cite{KK1}}, by which the method used in this article is suggested.
Besides, \eqref{Eq:qP6_JS} was derived from $q$-UC hierarchy in {\rm\cite{TM}}.
\end{remark}

The $q$-hypergeometric function ${}_n\phi_{n-1}$ is defined by the power series
\begin{equation*}
	{}_n\phi_{n-1}\left[\begin{array}{c}\alpha_1,\ldots,\alpha_{n-1},\alpha_n\\\beta_1,\ldots,\beta_{n-1}\end{array};q,t\right] = \sum_{k=0}^{\infty}\frac{(\alpha_1;q)_k\ldots(\alpha_{n-1};q)_k(\alpha_n;q)_k}{(\beta_1;q)_k\ldots(\beta_{n-1};q)_k(q;q)_k}t^k,
\end{equation*}
where $(\alpha;q)_k$ stands for the $q$-shifted factorial
\begin{equation*}
	(\alpha;q)_0 = 1,\quad (\alpha;q)_k = (1-\alpha)(1-q\alpha)\ldots(1-q^{k-1}\alpha)\quad (k\geq1).
\end{equation*}
Note that it reduces to the generalized hypergeometric function ${}_nF_{n-1}$ via a continuous limit $q\to1$.
We see that $x(t)={}_n\phi_{n-1}$ satisfies a $n$-th order linear $q$-difference equation
\begin{equation*}\begin{split}
	&\Bigl((1-q^{-1}\beta_1T_{q,t})\ldots(1-q^{-1}\beta_{n-1}T_{q,t})(1-T_{q,t})\\
	&\quad -t(1-\alpha_1T_{q,t})\ldots(1-\alpha_{n-1}T_{q,t})(1-\alpha_nT_{q,t})\Bigr)x(t) = 0,
\end{split}\end{equation*}
where $T_{q,t}$ stands for a $q$-shift operator such that $T_{q,t}x(t)=x(qt)$.

This article is organized as follows.
In Section \ref{Sec:qDS}, we formulate a $q$-DS hierarchy and its similarity reduction corresponding to the partition $(n,\ldots,n)$ of $mn$.
In Section \ref{Sec:qPainleve}, an explicit formula of $q$-$P_{(n,n)}$ is presented.
We also discuss a group of symmetries and a Lax form for $q$-$P_{(n,n)}$.
In Section \ref{Sec:qP6}, the $q$-Painlev\'{e} VI equation is derived from $q$-$P_{(2,2)}$.
In Section \ref{Sec:qHG}, we show that $q$-$P_{(n,n)}$ admits a particular solution in terms of ${}_n\phi_{n-1}$.

\section{$q$-DS hierarchy}\label{Sec:qDS}

In this section, we formulate a $q$-DS hierarchy and its similarity reduction corresponding to the partition $(n,\ldots,n)$ of $mn$.

Let $\lambda$ and $T_{\lambda}$ be a grading parameter and a corresponding $q$-shift operator such that $T_{\lambda}(\lambda)=q\lambda$.
We define $mn\times mn$ matrices $e_j,f_j,h_j$ $(j=1,\ldots,mn)$ by
\begin{equation*}\begin{array}{llll}
	e_j = E_{j,j+1}& f_j = E_{j+1,j}& h_j = E_{j,j}& (j\notin m\mathbb{Z}),\\[4pt]
	e_j = \lambda E_{j,j+1}& f_j = \lambda^{-1}E_{j+1,j}& h_j = E_{j,j}& (j\in m\mathbb{Z}; j\neq mn),\\[4pt]
	e_{mn}= \lambda E_{mn,1}& f_{mn} = \lambda^{-1}E_{1,mn}& h_{mn} = E_{mn,mn},
\end{array}\end{equation*}
where $E_{i,j}$ stands for the matrix with 1 on the $(i,j)$-th entry and zeros elsewhere.
We also set
\begin{equation*}
	e_{j,k} = e_je_{j+1}\ldots e_{j+k-1},\quad f_{j,k} = f_{j+k-1}\ldots f_{j+1}f_j,
\end{equation*}
where $e_{j+mn}=e_j$ and $f_{j+mn}=f_j$.
Note that
\begin{equation*}
	T_{\lambda}(e_{j,m}) = qe_{j,m},\quad T_{\lambda}(f_{j,m}) = q^{-1}f_{j,m}\quad (j=1,\ldots,mn).
\end{equation*}

Let us consider a set of matrices $\{\Lambda_1,\ldots,\Lambda_m\}$ such that
\begin{equation*}
	\Lambda_i\Lambda_j = \Lambda_j\Lambda_i = O,\quad T_{\lambda}(\Lambda_i) = q\Lambda_i\quad (i,j=1,\ldots,m; i\neq j).
\end{equation*}
Explicitly, we can take
\begin{equation*}
	\Lambda_i = \sum_{j=1}^{n}e_{i+(j-1)m,m}\quad (i=1,\ldots,m).
\end{equation*}
Note that such choice of matrices is suggested by the previous work \cite{BK,KL,S2}.
In the following, we formulate a $q$-DS hierarchy by using such matrices.

Let $t_i$ and $T_i$ $(i=1,\ldots,m)$ be independent variables and corresponding $q$-shift operators such that
\begin{equation*}
	T_i(t_i) = qt_i,\quad T_i(t_j) = t_j\quad (j\neq i).
\end{equation*}
Under the condition $|q|>1$, we consider a $q$-Sato equation
\begin{equation}\label{Eq:Sato}
	T_i(Z)Z^{-1} = T_i(W)(I-\varepsilon t_i\Lambda_i)W^{-1}\quad (i=1,\ldots,m),
\end{equation}
with matrices of dependent variables
\begin{equation*}
	Z = \sum_{j=1}^{mn}z_{j,0}h_j + \sum_{j=1}^{mn}\sum_{k=1}^{\infty}z_{j,k}e_{j,k},\quad W = I + \sum_{j=1}^{mn}\sum_{k=1}^{\infty}w_{j,k}f_{j,k},
\end{equation*}
where $I$ stands for the identity matrix and $\varepsilon=1-q$.
We assume that
\begin{equation}\label{Eq:Sato_tr}
	z_{1,0}z_{2,0}\ldots z_{mn,0} = 1.
\end{equation}
Under the system \eqref{Eq:Sato}, we set
\begin{equation*}
	\Psi = W\prod_{i=1}^{m}\prod_{k=0}^{\infty}(I-q^{-k-1}\varepsilon t_i\Lambda_i).
\end{equation*}
Then we obtain a system of linear $q$-difference equations called a Lax form
\begin{equation}\label{Eq:Lax_DS}
	T_i(\Psi) = B_i\Psi\quad (i=1,\ldots,m),
\end{equation}
where
\begin{equation*}
	B_i = T_i(W)(I-\varepsilon t_i\Lambda_i)W^{-1}.
\end{equation*}
The compatibility condition of the system \eqref{Eq:Lax_DS} is equivalent to
\begin{equation}\label{Eq:DS}
	T_i(B_j)B_i = T_j(B_i)B_j\quad (i,j=1,\ldots,m).
\end{equation}
We call a system of $q$-difference equations \eqref{Eq:DS} a $q$-DS hierarchy.

Next, we formulate a similarity reduction of the $q$-DS hierarchy.
Under the system \eqref{Eq:DS}, we consider the following equation:
\begin{equation}\label{Eq:Sato_SR}
	T_{\lambda}(W) = q^{\rho}T_{1,\ldots,m}(W)q^{-\rho},\quad \rho = \sum_{i=1}^{m}\sum_{j=1}^{n}\rho_ih_{i+(j-1)m},
\end{equation}
where $T_{i,\ldots,m}=T_iT_{i+1}\ldots T_m$ and $\rho_1+\ldots+\rho_m=0$.
Note that
\begin{equation*}
	\rho\Lambda_i = \Lambda_i\rho\quad (i=1,\ldots,m).
\end{equation*}
We set
\begin{equation*}
	\Psi = W\left(\prod_{i=1}^{m}\prod_{k=0}^{\infty}(I-q^{-k-1}\varepsilon t_i\Lambda_i)\right)\lambda^{\rho}.
\end{equation*}
Then we obtain a Lax form
\begin{equation}\label{Eq:Lax_SR}
	T_{\lambda}(\Psi) = M\Psi,\quad T_i(\Psi) = B_i\Psi\quad (i=1,\ldots,m),
\end{equation}
where
\begin{equation*}
	M = q^{\rho}T_{2,\ldots,m}(B_1)T_{3,\ldots,m}(B_2)\ldots T_m(B_{m-1})B_m.
\end{equation*}
Note that the matrix $M$ is also given by
\begin{equation*}
	M = T_{\lambda}(W)q^{\rho}(1-\varepsilon t_1\Lambda_1)(1-\varepsilon t_2\Lambda_2)\ldots(1-\varepsilon t_m\Lambda_m)W^{-1}.
\end{equation*}
The compatibility condition of the system \eqref{Eq:Lax_SR} is equivalent to
\begin{equation}\label{Eq:DS_SR}
	T_i(M)B_i = T_{\lambda}(B_i)M,\quad T_i(B_j)B_i = T_j(B_i)B_j\quad (i,j=1,\ldots,m).
\end{equation}
We call a system of $q$-difference equations \eqref{Eq:DS_SR} a similarity reduction of the $q$-DS hierarchy.

\begin{remark}
Replacing $M\to I-\varepsilon M$ and $B_i\to I-\varepsilon t_iB_i$, we can rewrite the system \eqref{Eq:DS_SR} into
\begin{equation*}
	[D_{\lambda}-M,D_i-B_i] = 0,\quad [D_i-B_i,D_j-B_j] = 0\quad (i,j=1,\ldots,m),
\end{equation*}
where
\begin{equation*}
	D_{\lambda} = \frac{1-T_{\lambda}}{\varepsilon\lambda},\quad D_i = \frac{1-T_i}{\varepsilon t_i}.
\end{equation*}
Via a continuous limit $q\to1$, it reduces to the similarity reduction of the DS hierarchy of type $A^{(1)}_{mn-1}$ given in {\rm\cite{S2}}.
\end{remark}

\section{Higher order $q$-Painlev\'{e} system}\label{Sec:qPainleve}

In this section, we present an explicit formula of the similarity reduction \eqref{Eq:DS_SR} corresponding to the partition $(n,n)$, where $n\geq2$, in terms of dependent variables.
We also discuss a group of symmetries and a Lax form for $q$-$P_{(n,n)}$.

Consider a Lax form
\begin{equation}\label{Eq:Lax_SR_2n}
	T_{\lambda}(\Psi) = M\Psi,\quad T_i(\Psi) = B_i\Psi\quad (i=1,2),
\end{equation}
with $2n\times2n$ matrices
\begin{equation*}\begin{split}
	M &= T_{\lambda}(W)q^{\rho}(1-\varepsilon t_1\Lambda_1)(1-\varepsilon t_2\Lambda_2)W^{-1},\\
	B_i &= T_i(W)(1-\varepsilon t_i\Lambda_i)W^{-1}\quad (i=1,2),
\end{split}\end{equation*}
where
\begin{equation*}
	\Lambda_1 = \sum_{j=1}^{n}e_{2j-1,2},\quad \Lambda_2 = \sum_{j=1}^{n}e_{2j,2},\quad \rho = \sum_{j=1}^{n}\rho_1(h_{2j-1}-h_{2j}),
\end{equation*}
and
\begin{equation*}
	W = I + \sum_{j=1}^{2n}\sum_{k=1}^{\infty}w_{j,k}f_{j,k}.
\end{equation*}
Then the compatibility condition of the system \eqref{Eq:Lax_SR_2n} is equivalent to the similarity reduction
\begin{equation}\label{Eq:DS_SR_2n}
	T_1(B_2)B_1 = T_2(B_1)B_2,\quad T_i(M)B_i = T_{\lambda}(B_i)M\quad (i=1,2).
\end{equation}
The matrices $B_1$, $B_2$ and $M$ can be expressed in terms of dependent variables $w_j=w_{j,1}$ and parameters $\kappa_j$ $(j=1,\ldots,2n)$.
We assume that indices of $w_j$ are congruent modulo $2n$, namely $w_j=w_{j+2n}$.
The equation \eqref{Eq:Sato_SR} implies
\begin{equation*}
	T_{1,2}(w_{2j-1}) = q^{2\rho_1}w_{2j-1},\quad T_{1,2}(w_{2j}) = q^{-2\rho_1-1}w_{2j}\quad (j=1,\ldots,n),
\end{equation*}
from which we obtain
\begin{equation*}\begin{split}
	M &= \sum_{j=1}^{2n}(1-\varepsilon\kappa_j)h_j + \sum_{j=1}^{n}(q^{\rho_1}\varepsilon t_1w_{2j}-q^{-\rho_1-1}\varepsilon t_2w_{2j-2})e_{2j-1}\\
	&\quad + \sum_{j=1}^{n}(q^{-\rho_1}\varepsilon t_2w_{2j+1}-q^{\rho_1}\varepsilon t_1w_{2j-1})e_{2j} - q^{\rho_1}\varepsilon t_1\Lambda_1 - q^{-\rho_1}\varepsilon t_2\Lambda_2,
\end{split}\end{equation*}
and
\begin{equation*}
	B_i = \sum_{j=1}^{2n}u_{i,j}h_j + \sum_{j=1}^{2n}v_{i,j}e_j - \varepsilon t_i\Lambda_i\quad (i=1,2),
\end{equation*}
where
\begin{equation*}\begin{split}
	&u_{1,2j-1} = \frac{q^{-\rho_1}(1-\varepsilon\kappa_{2j-1})}{1+q^{-2\rho_1-1}\varepsilon t_2w_{2j-2}T_1(w_{2j-1})},\quad
	u_{1,2j} = 1 + \varepsilon t_1T_1(w_{2j-1})w_{2j},\\
	&v_{1,2j-1} = \varepsilon t_1w_{2j},\quad
	v_{1,2j} = -\varepsilon t_1T_1(w_{2j-1}),\\
	&u_{2,2j-1} = 1 + \varepsilon t_2T_2(w_{2j-2})w_{2j-1},\quad
	u_{2,2j} = \frac{q^{\rho_1}(1-\varepsilon\kappa_{2j})}{1+q^{2\rho_1}\varepsilon t_1w_{2j-1}T_2(w_{2j})},\\
	&v_{2,2j-1} = -\varepsilon t_2T_2(w_{2j-2}),\quad
	v_{2,2j} = \varepsilon t_2w_{2j+1}.
\end{split}\end{equation*}
Note that the equation \eqref{Eq:Sato_tr} implies
\begin{equation*}
	\prod_{j=1}^{2n}(1-\varepsilon\kappa_j) = 1,\quad \prod_{j=1}^{2n}u_{i,j} = 1\quad (i=1,2).
\end{equation*}

We now assume that $t_2=1$.
We also consider a change of variables and parameters
\begin{equation*}\begin{split}
	&p=q^n,\quad t = q^{2n\rho_1}t_1^n,\quad z = q^{-n(4\rho_1-n+1)/2}(\epsilon\lambda)^n,\\
	&a_j = q^{-\rho_1+j-1}(1-\varepsilon\kappa_{2j-1}),\quad b_j = q^{-\rho_1+j-1}(1-\varepsilon\kappa_{2j}),\\
	&x_j(t) = q^{-2(j-1)\rho_1}t_1^{-j+1}w_{2j-1},\quad y_j(t) = q^{(j-1)(2\rho_1+1)}t_1^j\varepsilon w_{2j},
\end{split}\end{equation*}
for $j=1,\ldots,n$.
Replacing $p\to q$, we can describe the similarity reduction \eqref{Eq:DS_SR_2n} as follows.
\begin{theorem}
The dependent variables $x_j(t)$, $y_j(t)$ $(j=1,\ldots,n)$ satisfy a system of $q$-difference equations
\begin{equation}\begin{split}\label{Eq:Painleve}
	x_j(t) - x_{j-1}(t) &= \frac{a_jx_j(qt)}{1+x_j(qt)y_{j-1}(t)} - \frac{b_{j-1}x_{j-1}(qt)}{1+x_{j-1}(qt)y_{j-1}(t)},\\
	y_j(qt) - y_{j-1}(qt) &= \frac{b_jy_j(t)}{1+x_j(qt)y_j(t)} - \frac{a_jy_{j-1}(t)}{1+x_j(qt)y_{j-1}(t)},
\end{split}\end{equation}
for $j=1,\ldots,n$, where
\begin{equation*}
	b_0 = q^{-1}b_n,\quad x_0(t) = tx_n(t),\quad y_0(t) = q^{-1}t^{-1}y_n(t),
\end{equation*}
with relations
\begin{equation*}\begin{split}
	\prod_{j=1}^{n}a_j\frac{1+x_j(qt)y_j(t)}{1+x_j(qt)y_{j-1}(t)} = q^{(n-1)/2}.
\end{split}\end{equation*}
\end{theorem}

We denote the $q$-Painlev\'{e} system \eqref{Eq:Painleve} by $q$-$P_{(n,n)}$.
As is seen in the next section, $q$-$P_{(2,2)}$ coincides with the $q$-Painlev\'{e} VI equation \eqref{Eq:qP6_JS}.
Note that $q$-$P_{(n,n)}$ reduces to the coupled Painlev\'{e} VI system given in \cite{S2,T4} via a continuous limit $q\to1$.

\begin{remark}
In the previous work {\rm\cite{T1}}, the higher order generalizations of the $q$-Painlev\'{e} VI equation was presented.
Its relationship to $q$-$P_{(n,n)}$ has not been clarified yet.
However, we conjecture that $q$-$P_{(n,n)}$ coincides with the $q$-Painlev\'{e} equation of {\rm\cite{T1}} with $(n,n)$ periodicity and $(|I|,|J|)=(2,2)$.
\end{remark}

The system $q$-$P_{(n,n)}$ admits the affine Weyl group symmetry of type $A_{2n-1}^{(1)}$.
We describe its action on the dependent variables and parameters.
Recall that an extended affine Weyl group $\widetilde{W}(A^{(1)}_{2n-1})$ is generated by the transformations $r_0,\ldots,r_{2n-1}$ and $\pi$ with the fundamental relations
\begin{equation*}\begin{split}
	&r_i^2=1,\quad
	(r_ir_j)^{2-a_{i,j}}=1\quad (i,j=0,\ldots,2n-1; i\neq j),\\
	&\pi^{2n} = 1,\quad
	\pi r_i = r_{i+1}\pi,\quad
	\pi r_{2n-1} = r_0\pi\quad (i=0,\ldots,2n-2),
\end{split}\end{equation*}
where
\begin{equation*}\begin{array}{llll}
	a_{i,i}=2& (i=0,\ldots,2n-1),\\[4pt]
	a_{i,i+1}=a_{2n-1,0}=a_{i+1,i}=a_{0,2n-1}=-1& (i=0,\ldots,2n-2),\\[4pt]
	a_{i,j}=0& (\text{otherwise}).
\end{array}\end{equation*}

\begin{theorem}
Let $r_j$ $(j=0,\ldots,2n-1)$ be birational transformations defined by
\begin{equation*}\begin{split}
	&r_{2j-2}(a_j) = b_{j-1},\quad r_{2j-2}(b_{j-1}) = a_j,\\
	&r_{2j-2}(x_{j-1}(t)) = x_{j-1}(t),\quad r_{2j-2}(y_{j-1}(y)) = y_{j-1}(t) + \frac{b_{j-1}-a_j}{x_j(t)-x_{j-1}(t)},\\
	&r_{2j-2}(a_i) = a_i,\quad r_{2j-2}(b_{i-1}) = b_{i-1},\\
	&r_{2j-2}(x_{i-1}(t)) = x_{i-1}(t),\quad r_{2j-2}(y_{i-1}(y)) = y_{i-1}(t)\quad (i\neq j),
\end{split}\end{equation*}
and
\begin{equation*}\begin{split}
	&r_{2j-1}(a_j) = b_j\quad r_{2j-1}(b_j) = a_j,\\
	&r_{2j-1}(x_j(t)) = x_j(t) + \frac{a_j-b_j}{y_j(t)-y_{j-1}(t)},\quad r_{2j-1}(y_j(t)) = y_j(t),\\
	&r_{2j-1}(a_i) = a_i,\quad r_{2j-1}(b_i) = b_i,\\
	&r_{2j-1}(x_i(t)) = x_i(t),\quad r_{2j-1}(y_i(y)) = y_i(t)\quad (i\neq j).
\end{split}\end{equation*}
Also let $\pi$ be a birational transformation defined by
\begin{equation*}\begin{split}
	&a_i \to b_i,\quad b_i \to a_{i+1},\quad a_n \to b_n,\quad b_n \to qa_1,\\
	&x_i(t) \to y_i(qt),\quad y_i(t) \to x_{i+1}(qt),\\
	&x_n(t) \to y_n(qt),\quad y_n(t) \to \frac{x_1(qt)}{qt},\quad t \to \frac{1}{q^2t}\quad (i\neq n).
\end{split}\end{equation*}
Then $q$-$P_{(n,n)}$ is invariant under actions of those transformations.
Furthermore, the group of symmetries $\langle r_0,\ldots,r_{2n+1},\pi\rangle$ is isomorphic to an extended affine Weyl group $\widetilde{W}(A_{2n+1}^{(1)})$.
\end{theorem}

We rewrite the Lax form \eqref{Eq:Lax_SR_2n} into a simpler one.
Consider a gauge transformation
\begin{equation*}
	\widetilde{\Psi}(z,t) = {\lambda}^{-\rho_1}\sum_{j=1}^{n}q^{(j-1)(j-2)/2}\left((\varepsilon t_1\lambda)^{j-1}h_{2j-1}+(q^{-2\rho_1}\varepsilon\lambda)^{j-1}h_{2j}\right)\Psi.
\end{equation*}
Then we obtain a Lax form for $q$-$P_{(n,n)}$
\begin{equation}\label{Eq:Lax_Painleve}
	\widetilde{\Psi}(qz,t) = \widetilde{M}(z,t)\widetilde{\Psi}(z,t),\quad \widetilde{\Psi}(z,qt) = \widetilde{B}(z,t)\widetilde{\Psi}(z,t),
\end{equation}
with $2n\times2n$ matrices
\begin{equation*}\begin{split}
	\widetilde{M}(z,t) &= \sum_{j=1}^{n}a_jE_{2j-1,2j-1} + \sum_{j=1}^{n}b_jE_{2j,2j} + (y_1(t)-q^{-1}t^{-1}y_n(t))E_{1,2}\\
	&\quad + \sum_{j=1}^{n-1}(x_{j+1}(t)-x_j(t))E_{2j,2j+1} + \sum_{j=2}^{n}(y_j(t)-y_{j-1}(t))E_{2j-1,2j}\\
	&\quad + (x_1(t)-tx_n(t))zE_{2n,1} - \sum_{j=1}^{2n-2}E_{j,j+2} - tzE_{2n-1,1} - zE_{2n,2},
\end{split}\end{equation*}
and
\begin{equation*}\begin{split}
	\widetilde{B}(z,t) &= \sum_{j=1}^{n}\frac{a_j}{1+x_j(qt)y_{j-1}(t)}E_{2j-1,2j-1} + \sum_{j=1}^{n}(1+x_j(qt)y_j(t))E_{2j,2j}\\
	&\quad + \sum_{j=1}^{n}y_j(t)E_{2j-1,2j} - \sum_{j=1}^{n-1}x_j(qt)E_{2j,2j+1} - tx_n(qt)zE_{2n,1}\\
	&\quad - \sum_{j=1}^{n-1}E_{2j-1,2j+1} - tzE_{2n-1,1},
\end{split}\end{equation*}
where $E_{i,j}$ stands for the matrix unit.
The group of symmetries defined above is derived from transformations for the Lax form \eqref{Eq:Lax_Painleve}
\begin{equation*}
	r_j(\widetilde{\Psi}(z,t)) = R_i(z,t)\widetilde{\Psi}(z,t)\quad (j=0,\ldots,2n-1),
\end{equation*}
where
\begin{equation*}\begin{split}
	R_0(z,t) &= I + \frac{b_n-qa_1}{x_1(t)-tx_n(t)}z^{-1}E_{1,2n},\\
	R_{2j-2}(z,t) &= I + \frac{b_{j-1}-a_j}{x_j(t)-x_{j-1}(t)}E_{2j-1,2j-2}\quad (j=2,\ldots,n),\\
	R_{2j-1}(z,t) &= I + \frac{a_j-b_j}{y_j(t)-y_{j-1}(t)}E_{2j,2j-1}\quad (j=1,\ldots,n).
\end{split}\end{equation*}
Note that a construction of those transformations is suggested by the previous works \cite{KNY1,NY2}.

\section{$q$-Painlev\'{e} VI equation}\label{Sec:qP6}

The system $q$-$P_{(2,2)}$ is given as the compatibility condition of the Lax form
\begin{equation}\label{Eq:Lax_Painleve_2_4}
	\Psi_4(qz,t) = M_4(z,t)\Psi_4(z,t),\quad \Psi_4(z,qt) = B_4(z,t)\Psi_4(z,t),
\end{equation}
where
\begin{equation*}\begin{split}
	M_4(z,t) &= \begin{bmatrix}a_1&y_1(t)-\frac{y_2(t)}{qt}&-1&0\\0&b_1&x_2(t)-x_1(t)&-1\\-tz&0&a_2&y_2(t)-y_1(t)\\\{x_1(t)-tx_2(t)\}z&-z&0&b_2\end{bmatrix},\\
	B_4(z,t) &= \begin{bmatrix}\frac{qta_1}{qt+x_1(qt)y_2(t)}&y_1(t)&-1&0\\0&1+x_1(qt)y_1(t)&-x_1(qt)&0\\-tz&0&\frac{a_2}{1+x_2(qt)y_1(t)}&y_2(t)\\-tx_2(qt)z&0&0&1+x_2(qt)y_2(t)\end{bmatrix}.
\end{split}\end{equation*}
In this section, we derive the system \eqref{Eq:Lax_JS} from the Lax form \eqref{Eq:Lax_Painleve_2_4} with the aid of a $q$-Laplace transformation (cf. \cite{H,KK1}).

\begin{remark}
A $4\times4$ Lax form of the $q$-Painlev\'{e} VI equation similar to \eqref{Eq:Lax_Painleve_2_4} was already presented in {\rm\cite{PNGR}}.
\end{remark}

First, we consider a gauge transformation
\begin{equation*}
	\Psi_4^*(z,t) = \tau_1(z,t)\Psi_4(z,t),
\end{equation*}
where $\tau_1(z,t)$ is a function such that
\begin{equation*}
	\frac{\tau_1(qz,t)}{\tau_1(z,t)} = \frac{1}{qa_1},\quad \frac{\tau_1(z,qt)}{\tau_1(z,t)} = \frac{qt+x_1(qt)y_2(t)}{qta_1}.
\end{equation*}
Then the Lax form \eqref{Eq:Lax_Painleve_2_4} is transformed into
\begin{equation}\label{Eq:Lax_Painleve_2_4_scaling}
	\Psi_4^*(qz,t) = M_4^*(z,t)\Psi_4^*(z,t),\quad \Psi_4^*(z,qt) = B_4^*(z,t)\Psi_4^*(z,t),
\end{equation}
where
\begin{equation*}
	M_4^*(z,t) = \frac{1}{qa_1}M_4(z,t),\quad B_4^*(z,t) = \frac{qt+x_1(qt)y_2(t)}{qta_1}B_4(z,t).
\end{equation*}
We set
\begin{equation*}
	M_4^*(z,t) = M_{4,0}^*(t) + zM_{4,1}^*(t),\quad B_4^*(z,t) = B_{4,0}^*(t) + zB_{4,1}^*(t).
\end{equation*}
Next, we consider a $q$-Laplace transformation
\begin{equation*}
	z\Psi_4^*(z) \to \frac{\Phi_4(\zeta)-\Phi_4(q^{-1}\zeta)}{\varepsilon\zeta},\quad \Psi_4^*(qz) \to q^{-1}\Phi_4(q^{-1}\zeta),
\end{equation*}
where $\varepsilon=1-q$.
Then the Lax form \eqref{Eq:Lax_Painleve_2_4_scaling} is transformed into
\begin{equation}\label{Eq:Lax_Painleve_2_4_Laplace}
	\Phi_4(q^{-1}\zeta,t) = N_4(\zeta,t)\Phi_4(\zeta,t),\quad
	\Phi_4(\zeta,qt) = C_4(\zeta,t)\Phi_4(z,t).
\end{equation}
where
\begin{equation*}\begin{split}
	N_4(\zeta,t) &= \left(q^{-1}\varepsilon\zeta I+M_{4,1}^*(t)\right)^{-1}\left(\varepsilon\zeta M_{4,0}^*(t)+M_{4,1}^*(t)\right),\\
	C_4(\zeta,t) &= B_{4,0}^*(t) + \varepsilon^{-1}\zeta^{-1}B_{4,1}^*(t)\left(I-N_4(\zeta,t)\right).
\end{split}\end{equation*}
Denoting $\zeta^{-1}$ by $z$, we can restrict the Lax form \eqref{Eq:Lax_Painleve_2_4_Laplace} to the one with $3\times3$ matrices
\begin{equation}\label{Eq:Lax_Painleve_2_3}
	\Psi_3(qz,t) = M_3(z,t)\Psi_3(z,t),\quad \Psi_3(z,qt) = B_3(z,t)\Psi_3(z,t),
\end{equation}
thanks to the following lemma.
\begin{lemma}
For each of the matrices $N_4(\zeta,t)$ and $C_4(\zeta,t)$, the first column equals to the fundamental vector ${}^t[1,0,0,0]$.
\end{lemma}

\begin{remark}
The Lax form \eqref{Eq:Lax_Painleve_2_3} coincides with the one for the $q$-$\mathfrak{gl}_3$ hierarchy given in {\rm\cite{KK1}}.
\end{remark}

In a similar way, we can reduce the Lax form \eqref{Eq:Lax_Painleve_2_3} to the one with $2\times2$ matrices.
We consider a gauge transformation
\begin{equation*}
	\Psi_3^*(z,t) = \tau_2(z,t)\Psi_3(z,t),
\end{equation*}
where $\tau_2(z,t)$ is a function such that
\begin{equation*}
	\frac{\tau_2(qz,t)}{\tau_2(z,t)} = \frac{a_1}{qb_1},\quad \frac{\tau_2(z,qt)}{\tau_2(z,t)} = \frac{qta_1}{(qt+x_1(qt)y_2(t))(1+x_1(qt)y_1(t))}.
\end{equation*}
We also consider a $q$-Laplace transformation
\begin{equation*}
	z\Psi_3^*(z) \to \frac{\Phi_3(\zeta)-\Phi_3(q^{-1}\zeta)}{\varepsilon\zeta},\quad \Psi_3^*(qz) \to q^{-1}\Phi_3(q^{-1}\zeta).
\end{equation*}
Then the Lax form \eqref{Eq:Lax_Painleve_2_3} is transformed into
\begin{equation}\label{Eq:Lax_Painleve_2_3_Laplace}
	\Phi_3(q^{-1}\zeta,t) = N_3(\zeta,t)\Phi_3(\zeta,t),\quad
	\Phi_3(\zeta,qt) = C_3(\zeta,t)\Phi_3(z,t).
\end{equation}
Denoting $\varepsilon^2a_1b_1\zeta$ by $z$, we can restrict the Lax form \eqref{Eq:Lax_Painleve_2_3_Laplace} to the one with $2\times2$ matrices
\begin{equation}\label{Eq:Lax_Painleve_2_2}
	\Psi_2(q^{-1}z,t) = M_2(z,t)\Psi_2(z,t),\quad \Psi_2(z,qt) = B_2(z,t)\Psi_2(z,t),
\end{equation}
thanks to the following lemma.
\begin{lemma}
For each of the matrices $N_3(\zeta,t)$ and $C_3(\zeta,t)$, the first column equals to the fundamental vector ${}^t[1,0,0]$.
\end{lemma}

The matrices $M_2(z,t)$ and $B_2(z,t)$ are of the form
\begin{equation*}\begin{split}
	M_2(z,t) = \frac{M_{2,0}(t)+zM_{2,1}(t)+z^2M_{2,2}(t)}{(z-t)(z-1)},\quad B_2(z,t) = \frac{B_{2,0}(t)+zB_{2,1}(t)}{z-t},
\end{split}\end{equation*}
where
\begin{equation*}\begin{split}
	M_{2,2}(t) &= \frac{1}{b_1}\begin{bmatrix}a_2&y_2(t)-y_1(t)\\0&b_2\end{bmatrix},\\
	B_{2,1}(t) &= \frac{1}{1+x_1(qt)y_1(t)}\begin{bmatrix}\frac{a_2}{1+x_2(qt)y_1(t)}&y_2(t)\\0&1+x_2(qt)y_2(t)\end{bmatrix}.
\end{split}\end{equation*}
In order to derive the system \eqref{Eq:Lax_JS}, we consider a gauge transformation
\begin{equation*}
	Y(z,t) = \frac{(q^{-1}tz^{-1};q^{-1})_{\infty}(qz;q)_{\infty}^2}{(q^{-1}z^{-1};q^{-1})_{\infty}}\begin{bmatrix}\tau_3(t)&0\\0&\tau_4(t)\end{bmatrix}\begin{bmatrix}1&\frac{y_2(t)-y_1(t)}{a_2-b_2}\\0&1\end{bmatrix}\Psi_2(z,t).
\end{equation*}
where $\tau_3(z,t)$ and $\tau_4(z,t)$ are functions such that
\begin{equation*}\begin{split}
	&\frac{\tau_3(q^{-1}z,t)}{\tau_3(z,t)} = b_1,\quad \frac{\tau_3(z,qt)}{\tau_3(z,t)} = \frac{(1+x_1(qt)y_1(t))(1+x_2(qt)y_1(t))}{a_2},\\
	&\frac{\tau_4(q^{-1}z,t)}{\tau_4(z,t)} = b_1,\quad \frac{\tau_4(z,qt)}{\tau_4(z,t)} = \frac{1+x_1(qt)y_1(t)}{1+x_2(qt)y_2(t)}.
\end{split}\end{equation*}
Then the Lax form \eqref{Eq:Lax_Painleve_2_2} is transformed into
\begin{equation}\label{Eq:Lax_qP6}
	Y(q^{-1}z,t) = \widetilde{M}_2(z,t)Y(z,t),\quad Y(z,t) = \frac{\widetilde{B}_2(z,q^{-1}t)}{z}Y(z,q^{-1}t).
\end{equation}
Here the coefficient matrices satisfy
\begin{equation*}\begin{split}
	&\widetilde{M}_2(z,t) = \widetilde{M}_{2,0}(t) + z\widetilde{M}_{2,1}(t) + z^2\widetilde{M}_{2,2},\\
	&\widetilde{M}_{2,2}=\begin{bmatrix}a_2&0\\0&b_2\end{bmatrix},\quad \text{$\widetilde{M}_{2,0}(t)$ has eigenvalues $ta_1,tb_1$},\\
	&\det\widetilde{M}_2(z,t) = a_2b_2(z-t)(z-a_1b_1q^{-1/2}t)(z-1)(z-a_2^{-1}b_2^{-1}q^{1/2}),
\end{split}\end{equation*}
and
\begin{equation*}\begin{split}
	&\widetilde{B}_2(z,q^{-1}t) = \widetilde{B}_{2,0}(q^{-1}t) + zI_2,\\
	&\det\widetilde{B}_2(z,q^{-1}t)=(z-q^{-1}t)(z-q^{-1}a_1b_1q^{-1/2}t).
\end{split}\end{equation*}
We can show them by direct computations; we do not state its detail here.
By setting
\begin{equation*}\begin{split}
	&\alpha_1 = 1,\quad \alpha_2 = \frac{a_1b_1}{q^{1/2}},\quad \alpha_3 = 1,\quad \alpha_4 = \frac{q^{1/2}}{a_2b_2},\\
	&\beta_1 = \frac{b_1}{q^{1/2}},\quad \beta_2 = \frac{a_1}{q^{1/2}},\quad \beta_3 = \frac{q}{a_2},\quad \beta_4 = \frac{1}{b_2},
\end{split}\end{equation*}
we arrive at
\begin{theorem}
The matrix-valued function $Y(z,t)$ solves the system \eqref{Eq:Lax_JS}.
\end{theorem}

\begin{corollary}
Under the system $q$-$P_{(2,2)}$, we set
\begin{equation*}
	x(t) = \frac{t(x_2(t)-x_1(t))\xi_1(t)}{\xi_2(t)},\quad y(t) = \frac{x_2(qt)(qt+x_1(qt)y_2(t))\psi_1(t)}{(1+x_2(qt)y_2(t))\psi_2(t)},
\end{equation*}
where
\begin{equation*}\begin{split}
	\xi_1(t) &= qtx_1(t)y_1(t) - x_1(t)y_2(t) - qtx_2(t)y_1(t) + x_2(t)y_2(t) - (b_1-a_1)qt,\\
	\xi_2(t) &= (tx_2(t)-x_1(t))(x_2(t)-x_1(t))(y_2(t)-qty_1(t))\\
	&\quad + (b_1-a_1)qtx_1(t) + \{(a_2-b_1)t-(a_2-a_1)\}qtx_2(t),\\
	\psi_1(t) &= (1-a_1b_1q^{1/2}t)x_2(qt)y_2(t) + qt - a_1b_1q^{1/2}t,\\
	\psi_2(t) &= a_2(1-a_1b_1q^{1/2}t)x_1(qt)x_2(qt)y_2(t)\\
	&\quad + a_1(qt-a_2b_1q^{1/2}t)x_1(qt) + (a_2-a_1)qtx_2(qt).
\end{split}\end{equation*}
Then those variables satisfy the $q$-Painlev\'{e} VI equation \eqref{Eq:qP6_JS}.
\end{corollary}

Note that  $q$-$P_{(2,2)}$ gives explicit formulas of $x_1(qt)$ and $x_2(qt)$ by
\begin{equation*}\begin{split}
	x_1(qt) &= \frac{-\xi_3(t)}{\xi_3(t)y_1(t)+(qt-a_1b_1q^{1/2}t)y_1(t)-(1-a_1b_1q^{1/2}t)y_2(t)},\\
	x_2(qt) &= \frac{-\xi_4(t)}{\xi_4(t)y_1(t)+a_2(qt-a_1b_1q^{1/2}t)y_1(t)-a_2(1-a_1b_1q^{1/2}t)y_2(t)},
\end{split}\end{equation*}
where
\begin{equation*}\begin{split}
	\xi_3(t) &= a_1q^{1/2}t(x_2(t)-x_1(t))(y_2(t)-y_1(t)) - qt + a_1a_2q^{1/2}t,\\
	\xi_4(t) &= (x_2(t)-x_1(t))(y_2(t)-qty_1(t)) - b_1(qt-a_1a_2q^{1/2}t).
\end{split}\end{equation*}

\section{$q$-Hypergeometric function ${}_n\phi_{n-1}$}\label{Sec:qHG}

In this section, we show that $q$-$P_{(n,n)}$ admits a particular solution in terms of the $q$-hypergeometric function ${}_n\phi_{n-1}$.

\begin{proposition}
Under the system $q$-$P_{(n,n)}$, we consider a specialization
\begin{equation*}
	y_j(t) = 0\quad (j=1,\ldots,n),\quad \prod_{j=1}^{n}a_j = q^{(n-1)/2}.
\end{equation*}
Then a vector of the variables $\mathbf{x}(t)={}^t[x_1(t),\ldots,x_n(t)]$ satisfies a system of linear $q$-difference equations
\begin{equation}\label{Eq:HGE}
	\mathbf{x}(q^{-1}t) = \left(A_0+\frac{A_1}{1-q^{-1}t}\right)\mathbf{x}(t),
\end{equation}
with $n\times n$ matrices
\begin{equation*}
	A_0 = \sum_{j=1}^{n}b_jE_{j,j} + \sum_{i=1}^{n}\sum_{j=i+1}^{n}(b_j-a_j)E_{i,j},\quad A_1 = \sum_{i=1}^{n}\sum_{j=1}^{n}(a_j-b_j)E_{i,j}.
\end{equation*}
\end{proposition}

We always assume that
\begin{equation*}
	a_j \notin \mathbb{Z},\quad a_i - a_j \notin \mathbb{Z}\quad (i,j=1,\ldots,n; i\neq j).
\end{equation*}
Note that $A_0$ is an upper triangular matrix and
\begin{equation*}
	A_0 + A_1 = \sum_{j=1}^{n}a_jE_{j,j} + \sum_{i=1}^{n}\sum_{j=1}^{i-1}(a_j-b_j)E_{i,j},
\end{equation*}
is a lower triangular matrix.

We consider a formal power series of $\mathbf{x}(t)$ at $t=0$
\begin{equation*}
	\mathbf{x}(t) = t^{\log_qa_1}\sum_{k=0}^{\infty}(q^{-1}t)^k\mathbf{x}_k.
\end{equation*}
Recall that $|q|>1$, namely $|q^{-1}|<1$.
Substituting it into the system \eqref{Eq:HGE}, we obtain
\begin{equation}\begin{split}\label{Eq:HGE_rf}
	(A_0+A_1-a_1I)\mathbf{x}_0 = \mathbf{0},\quad (A_0+A_1-a_1q^{-k}I)\mathbf{x}_k = A_0\mathbf{x}_{k-1}\quad (k\geq1).
\end{split}\end{equation}
The matrices $A_0+A_1-a_1I$ and $A_0+A_1-a_1q^{-k}I$ are of rank $n-1$ and $n$, respectively.
Hence the recurrence formula \eqref{Eq:HGE_rf} admits one parameter family of solutions.
Its explicit formula is given by
\begin{equation*}
	\mathbf{x}_k = \left(\frac{b_1\ldots b_n}{a_1\ldots a_n}\right)^k\begin{bmatrix}x_{k,1}\\\vdots\\x_{k,n}\end{bmatrix}\quad (k\geq0),
\end{equation*}
where
\begin{equation*}\begin{split}
	x_{0,1} = 1,\quad x_{0,j} = \prod_{i=1}^{j-1}\frac{b_i-a_1}{a_{i+1}-a_1}\quad (j=2,\ldots,n),
\end{split}\end{equation*}
and
\begin{equation*}\begin{split}
	x_{k,j} = \frac{(q\frac{a_1}{b_1};q)_k\ldots(q\frac{a_1}{b_{j-1}};q)_k(\frac{a_1}{b_j};q)_k\ldots(\frac{a_1}{b_{n-1}};q)_k(\frac{a_1}{b_n};q)_k}{(q\frac{a_1}{a_2};q)_k\ldots(q\frac{a_1}{a_j};q)_k(\frac{a_1}{a_{j+1}};q)_k\ldots(\frac{a_1}{a_n};q)_k(q;q)_k}x_{0,j},
\end{split}\end{equation*}
for $k\geq1$.
Then we arrive at
\begin{theorem}
The system \eqref{Eq:HGE} admits a solution
\begin{equation*}
	\mathbf{x}(t) = t^{\log_qa_1}\begin{bmatrix}c_1\varphi_1(t)\\\vdots\\c_n\varphi_n(t)\end{bmatrix},
\end{equation*}
where
\begin{equation*}\begin{split}
	c_j &= \prod_{i=1}^{j-1}\frac{b_i-a_1}{a_{i+1}-a_1},\\
	\varphi_j(t) &= {}_n\phi_{n-1}\left[\begin{array}{c}q\frac{a_1}{b_1},\ldots,q\frac{a_1}{b_{j-1}},\frac{a_1}{b_j},\ldots,q\frac{a_1}{b_{n-1}},\frac{a_1}{b_n}\\q\frac{a_1}{a_2},\ldots,q\frac{a_1}{a_j},\frac{a_1}{a_{j+1}},\ldots,\frac{a_1}{a_n}\end{array};q^{-1},\frac{b_1\ldots b_n}{a_1\ldots a_n}q^{-1}t\right].
\end{split}\end{equation*}
\end{theorem}

\section*{Acknowledgement}
The author would like to express his gratitude to Dr. Kenta Fuji for fruitful discussion.
The author is also grateful to Professors Saburo Kakei, Tetsuya Kikuchi, Masatoshi Noumi, Teruhisa Tsuda and Yasuhiko Yamada for helpful comments and advices.



\begin{thebibliography}{9}
\bibitem{BK}
	M. J. Bergvelt and A. P. E. ten Kroode,
	Partitions, vertex operator constructions and multi-component KP equations,
	Pacific J. Math. \textbf{171} (1995) 23-88.
\bibitem{DS}
	V. G. Drinfel'd and V. V. Sokolov,
	Lie algebras and equations of Korteweg-de Vries type,
	J. Sov. Math. \textbf{30} (1985) 1975-2036.
\bibitem{FHM}
	L. Feh\'{e}r, J. Harnad and I. Marshall,
	Generalized Drinfeld-Sokolov reductions and KdV type hierarchies,
	Comm. Math. Phys. \textbf{154} (1993) 181-214.
\bibitem{FN}
	H. Flaschka and Newell,
	Monodromy- and spectrum-preserving deformations I,
	Comm. Math. Phys. \textbf{76} (1980) 65-116.
\bibitem{FS1}
	K. Fuji and T. Suzuki,
	The sixth Painlev\'{e} equation arising from $D_4^{(1)}$ hierarchy,
	J. Phys. A: Math. Gen. \textbf{39} (2006) 12073-12082.
\bibitem{FS2}
	K. Fuji and T. Suzuki,
	Drinfeld-Sokolov hierarchies of type $A$ and fourth order Painlev\'{e} systems,
	Funkcial. Ekvac. \textbf{53} (2010) 143-167.
\bibitem{FS3}
	K. Fuji and T. Suzuki,
	Higher order Painlev\'{e} systems of type $A$, Drinfeld-Sokolov hierarchies and Fuchsian systems,
	RIMS Kokyuroku Bessatsu \textbf{B30} (2012) 181-208.
\bibitem{GHM}
	M. F. de Groot, T. J. Hollowood and J. L. Miramontes,
	Generalized Drinfeld-Sokolov hierarchies,
	Comm. Math. Phys. \textbf{145} (1992) 57-84.
\bibitem{H}
	W. Hahn,
	Beitr\"{a}ge zur theorie der heineschen reihen
	Math. Nachr. \textbf{2} (1949) 340-379.
\bibitem{JS}
	M. Jimbo and H. Sakai,
	A $q$-analog of the sixth Painlev\'{e} equation,
	Let. Math. Phys. \textbf{38} (1996) 145-154.
\bibitem{KK1}
	S. Kakei and T. Kikuchi,
	A $q$-analogue of $\mathfrak{gl}_3$ hierarchy and $q$-Painlev\'{e} VI,
	J. Phys. A: Math. Gen. \textbf{39} (2006) 12179-12190.
\bibitem{KK2}
	S. Kakei and T. Kikuchi,
	The sixth Painlev\'{e} equation as similarity reduction of $\widehat{\mathfrak{gl}}_3$ generalized Drinfel'd-Sokolov hierarchy,
	Lett. Math. Phys. \textbf{79} (2007) 221-234.
\bibitem{KL}
	F. ten Kroode and J. van de Leur,
	Bosonic and fermionic realizations of the affine algebra $\widehat{\mathfrak{gl}}_n$,
	Comm. Math. Phys. \textbf{137} (1991) 67-107.
\bibitem{KNY1}
	K. Kajiwara, M. Noumi and Y. Yamada,
	Discrete dynamical systems with $W(A_{m-1}^{(1)}\times A_{n-1}^{(1)})$ symmetry,
	Lett. Math. Phys. \textbf{60} (2002) 211-219.
\bibitem{KNY2}
	K. Kajiwara, M. Noumi and Y. Yamada,
	$q$-Painlev\'{e} systems arising from $q$-KP hierarchy,
	Lett. Math. Phys. \textbf{62} (2003) 259-268.
\bibitem{NY1}
	M. Noumi and Y. Yamada,
	Symmetries in the fourth Painlev\'{e} equation and Okamoto polynomials,
	Nagoya Math. J. \textbf{153} (1999) 53-86.
\bibitem{NY2}
	M. Noumi and Y. Yamada,
	Birational Weyl group action arising from a nilpotent Poisson algebra,
	in Physics and Combinatorics 1999, Proceedings of the Nagoya 1999 International Workshop, ed. A.N.Kirillov, A.Tsuchiya and H.Umemura, (World Scientific, 2001) 287-319.
\bibitem{PNGR}
	V. G. Papageorgiou, F. W. Nijhoff, B. Grammaticos and A. Ramani,
	Isomonodromic deformation problems for discrete analogues of Painlev\'{e} equations,
	Phys. Lett. A \textbf{164} (1992) 57-64.
\bibitem{S1}
	T. Suzuki,
	A particular solution of a Painlev\'{e} system in terms of the hypergeometric function ${}_{n+1}F_n$,
	SIGMA \textbf{6} (2010) 078.
\bibitem{S2}
	T. Suzuki,
	A class of higher order Painlev\'{e} systems arising from integrable hierarchies of type $A$,
	Contemp. Math. \textbf{593} (2013) 125-141.
\bibitem{T1}
	T. Tsuda,
	On an integrable system of $q$-difference equations satisfied by the universal characters: its Lax formalism and an application to $q$-Painlev\'{e} equations,
	Comm. Math. Phys. \textbf{293} (2010) 347-359.
\bibitem{T2}
	T. Tsuda,
	From KP/UC hierarchies to Painlev\'{e} equations,
	Internat. J. Math. \textbf{23} (2012) 1250010.
\bibitem{T3}
	T. Tsuda,
	Hypergeometric solution of a certain polynomial Hamiltonian system of isomonodromy type,
	Quart. J. Math. \textbf{63} (2012) 489-505.
\bibitem{T4}
	T. Tsuda,
	UC hierarchy and monodromy preserving deformation,
	J. reine angew. Math., in press.
\bibitem{TM}
	T. Tsuda and T. Masuda,
	$q$-Painlev\'{e} VI equation arising from $q$-UC hierarchy,
	Comm. Math. Phys. \textbf{262} (2006) 595-609.
\end{thebibliography}
\end{document}